\documentclass[journal]{IEEEtran}
\usepackage{cite}

\ifCLASSINFOpdf
\else
\fi
\usepackage{bm}   
\usepackage{amsmath}
\usepackage{amssymb}
\usepackage{indentfirst}  
\usepackage{graphicx}  
\usepackage{subfigure}
\usepackage{multirow} 
\usepackage{booktabs} 
\usepackage{cite}  
\usepackage{threeparttable}   
\usepackage{algorithm} 
\usepackage{algpseudocode}
\usepackage{arydshln}
\usepackage{color}
\usepackage{microtype}
\usepackage{cases}
\usepackage{enumerate}

\newtheorem{theorem}{Theorem}
\newtheorem{proposition}{Proposition}
\newtheorem{definition}{Definition}

\newtheorem{remark}{Remark}

\newcommand{\tabincell}[2]{\begin{tabular}{@{}#1@{}}#2\end{tabular}}

\algtext*{EndWhile}
\algtext*{EndIf}
\usepackage{multirow}

\begin{document}
%
\title{Radiality Constraints for Resilient Reconfiguration of Distribution Systems:
Formulation and Application to Microgrid Formation}
%
%
%

\author{Shunbo~Lei,~\IEEEmembership{Member,~IEEE,}
        Chen~Chen,~\IEEEmembership{Member,~IEEE,}
        Yue~Song,~\IEEEmembership{Member,~IEEE,}
        and~Yunhe~Hou,~\IEEEmembership{Senior~Member,~IEEE}
\thanks{The work of Y. Song and Y. Hou was supported in part by
the Research Grant Council of Hong Kong through
the Theme-based Research Scheme under Project No. T23-701/14-N.
The work of Y. Hou was also supported in part by the National
Natural Science Foundation
of China under Grant 51677160,
and in part by the
Research Grant Council of Hong Kong under Grant GRF17207818.}
\thanks{S. Lei is with the Department of Electrical Engineering and
Computer~Science, University of Michigan, Ann Arbor, MI 48109 USA;
C.~Chen~is~with~the Energy Systems Division, Argonne National Laboratory,
Argonne,~IL~60439 USA;
Y. Song and Y. Hou are with the Department of Electrical~and~Electronic
Engineering, The University of Hong Kong, Hong Kong; Y. Hou is also
with The University of Hong Kong Shenzhen Institute of Research and Innovation,
Shenzhen 518027 China
(e-mail: shunbol@umich.edu; morningchen@anl.gov; yuesong@eee.hku.hk;
yhhou@eee.hku.hk).}}
\maketitle

\begin{abstract}
Network reconfiguration is an effective strategy for different purposes
of distribution systems (DSs), e.g., resilience enhancement. In particular,
DS automation, distributed generation integration and microgrid (MG)
technology development, etc., are empowering much more flexible
reconfiguration and operation
of the system, e.g., DSs or MGs with flexible boundaries.
However, the formulation of DS reconfiguration-related~optimization problems
to include those new flexibilities is non-trivial, especially for the issue of
topology, which has to be radial. That is, the existing methods of formulating the
radiality constraints can cause under-utilization of DS flexibilities. Thus, in this
work, we propose a new method for radiality constraints formulation that fully
enables the topological and some other related flexibilities of DSs, so that the
reconfiguration-related optimization problems can have extended feasibility
and enhanced optimality. Graph-theoretic supports are provided to certify its
theoretical validity. As integer variables~are
involved,~we~also~analyze~the~issues~of~tightness~and~compactness.
The proposed radiality constraints are
specifically applied to post-disaster MG formation, which is involved in
many DS resilience-oriented service restoration and/or infrastructure recovery
problems. The resulting new MG formation model, which allows more flexible
merge and/or separation of the sub-grids, etc., establishes superiority over
the models in the literature. Demonstrative case studies are conducted on
two test systems.

\end{abstract}


%
\IEEEpeerreviewmaketitle

\section{Introduction}
%
%
%
%
\IEEEPARstart{D}{istribution}
system (DS) reconfiguration is an effective and multi-function
strategy \cite{Baran1989, Kening2016, ChenXin2016}. Its optimization thus has been extensively
studied.
As most DSs have to operate with a radial topology,
the mathematical formulation of radiality constraints
has been specifically investigated \cite{Ji1996,Ram2005,Lav2012,Jab2012,Ahm2015}
(see more details in Section \ref{tightness_issue}).
This issue is resolved for conventional DS reconfiguration.
Nevertheless, now DSs can adopt much more adaptive reconfiguration and operation
owing to
the added flexibilities of automation equipment, distributed generations (DGs),
and microgrid (MG) components, etc.
In particular, regarding the topology issue, DSs and MGs now can have flexible
boundaries \cite{Arefifar2015, Kim2018, Manshadi2018}.
For example, the DS is split~into a
to-be-optimized number of MGs in~\cite{Manshadi2018}.
Such added~flexibilities are actually empowering more~resilient~reconfiguration~of~DSs.

However, existing methods to formulate radiality constraints
cannot fully include those new flexibilities
in optimizing~DS reconfiguration \cite{SSMA2018}.
That is,
the feasible region given by these formulations only refers to a subset of the
actual one.~DS~flexibilities thus will be underutilized and
less coordinated,
which is especially adverse for resilience-oriented~reconfiguration~\cite{Bie2017}.
For example, post-disaster reconfiguration
is faced with quite limited flexibilities due to many faults
caused by an extreme weather event~\cite{Lei2016}, etc.
In such cases, the faults also
complicates decision-making and even requires co-optimization~with
other recovery efforts~\cite{Ari2017}, etc.
Actually,~to~co-optimize~with the repairing sequence of the damaged parts
involves reconfiguring a DS network with a physical structure
evolving with the repairing variables. Existing methods to formulate radiality constraints
cannot properly handle
these situations.

This work proposes a new formulation of
radiality constraints that
fully enables
topological and some other related flexibilities
in DS reconfiguration-related optimization problems.~It~is~superior to
the literature's other attempts
with the same or similar aims \cite{SSMA2018} \cite{Arif2018}
(see comparisons in Section \ref{proposed_method}).
Graph-theoretic justifications are provided to
affirm the analytical validity of it.
As integer variables are involved,
the tightness and compactness issues are also analyzed \cite{UCF2013}.
Generally, adopting the proposed radiality constraints
in DS reconfiguration optimization problems can attain
extended feasibility and enhanced optimality.

\begin{figure}[t!]
  \begin{minipage}[b!]{0.18\textwidth}
    \centering
    \includegraphics[width=0.925in]{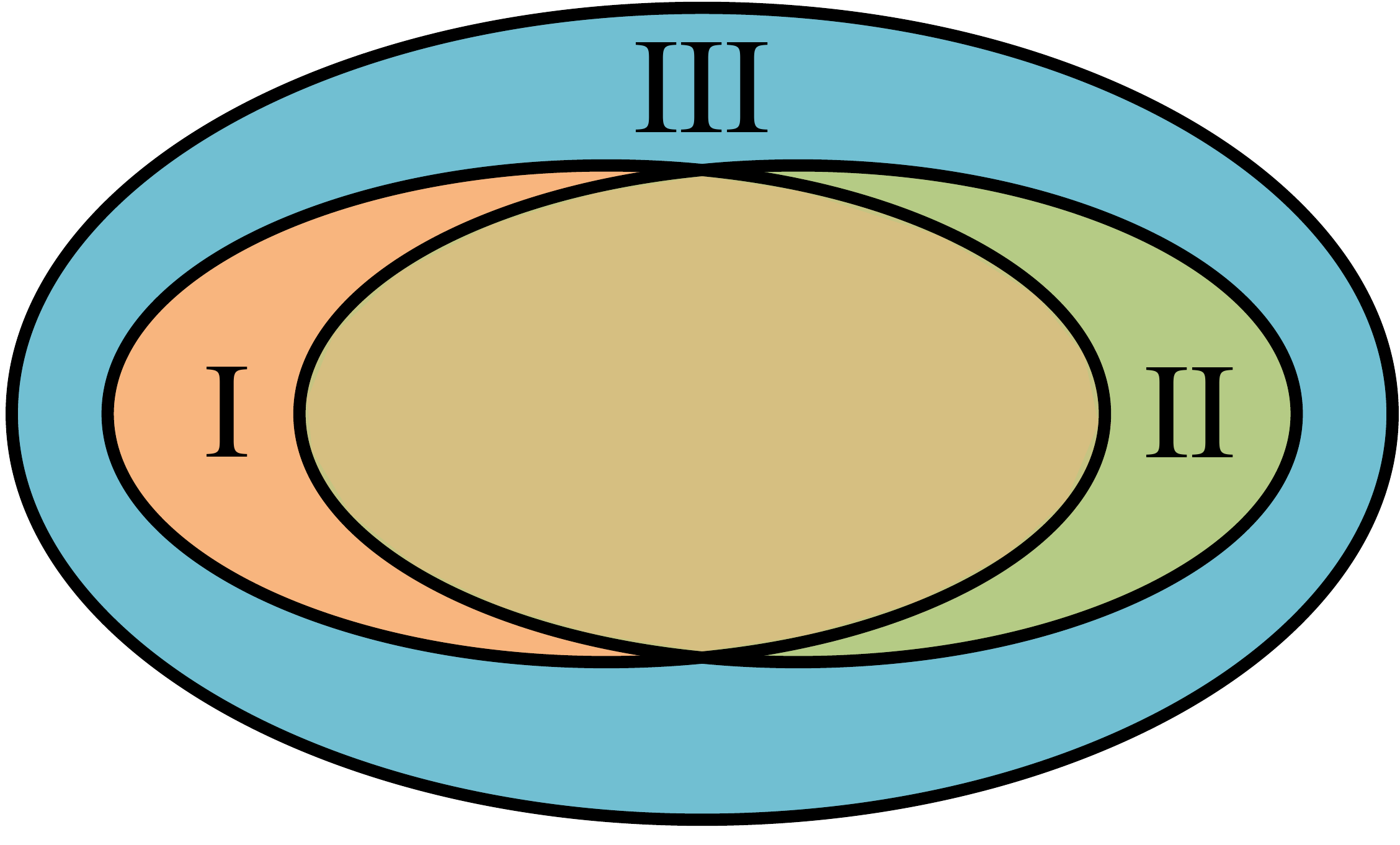}
    \label{feasibleRegions}
  \end{minipage}%
  \begin{minipage}[t!]{0.305\textwidth}
    \centering
    \setlength{\tabcolsep}{1.5pt}
    \begin{tabular}{rl}
      \small I: & \small Search space of the model in \cite{Chen2016} \\
      \small II:  & \small Search space of the model in \cite{Ding2017} \\
      \small III:    & \small Search space of the proposed model \\
    \end{tabular}
      \end{minipage}
  \caption{The search space (i.e., considered
  DS flexibilities) of our MG
  formation model and the models in \cite{Chen2016} \cite{Ding2017}.
  (Note: All comparisons in this work assume the same formulation
  other than the topology modeling. For example,
  loads are dispatchable in \cite{Ding2017}. We alter them to
  be non-dispatchable as in \cite{Chen2016}~and~here.)}
   \label{DifferentFeasibleRegions}
   \vspace{-1.9mm}
\end{figure}

For verification, the proposed radiality constraints
are~applied
to
construct a new optimization model for post-disaster MG~formation,
which reconfigures the DS to form MGs~energized~by DGs and/or other power sources.
It is an~essential~strategy~for many resilience-oriented DS restoration and/or
recovery~problems \cite{Ari2017, Arif2018, Lei2016, SSMA2018, Bie2017}.
Our proposed model again establishes superiority over the literature's
two groups of MG
formation models \cite{Chen2016} \cite{Ding2017}.
As compared in Fig. \ref{DifferentFeasibleRegions},
while their models exclude some DS flexibilities,
our model allows more flexible merge and separation of sub-grids, etc.
(See more details in Section~\ref{applicationMG}.)

In the following, Section \ref{proposed_method}
details the proposed method~of formulating radiality constraints.
Section
\ref{tightness_issue} analyzes its~tightness and compactness issues.
Sections \ref{applicationMG} and \ref{cases_MG}
apply~it~to~resilient MG formation.
Section \ref{conclusions_final}
provides the conclusion.

\section{Proposed Radiality Constraints}\label{proposed_method}

This work proposes to construct radiality constraints based~on two simple
graph-theoretic concepts
and their relationships.

First, the definition of \emph{spanning tree}, which may be already well-known, is still given
here to make this paper self-contained:

\begin{definition}\label{ST_defi}
  A spanning tree is a graph that connects all~the vertices and contains no cycles \cite{Bondy1976}.
\end{definition}

Second, \emph{spanning forest}, the other involved concept which is less common,
is defined as follows:

\begin{definition}\label{SF_defi}
  A spanning forest is a graph with no cycles~\cite{Bondy1976}.
\end{definition}

\begin{figure}[t!] 
  \centering
  \includegraphics[width=3.4in]{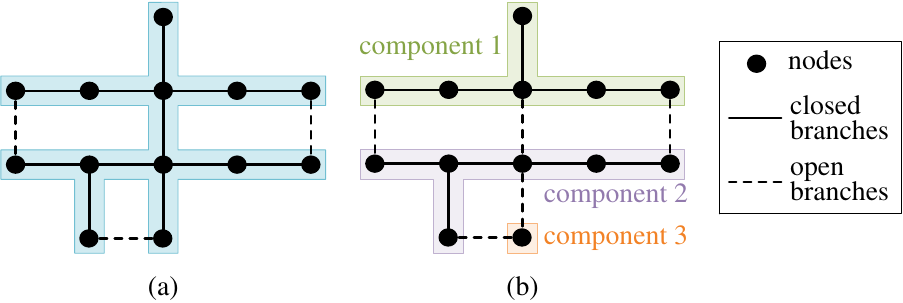}
  \vspace{-6pt}
  \caption{(a) A spanning tree; (b) A spanning forest. (Illustrated on the
  modified IEEE 13-node test system \cite{13system}.)}
  \vspace{-7pt}
  \label{STandSF}
\end{figure}

As in Fig.~\ref{STandSF}, a spanning forest is a graph whose \emph{connected components}
are spanning trees. Actually, a spanning forest~with $\kappa$ components can be called
as a $\kappa$-\emph{tree}. For example, Fig.~\ref{STandSF}(b) is a $3$-tree.
And, a $1$-tree is a spanning tree, e.g., Fig.~\ref{STandSF}(a).
Obviously, both spanning trees and spanning forests are radial.

Let $\kappa_1$ be the number of substations in the DS.
In light of Definitions 1 and 2, we can observe that normal DS
reconfiguration for loss reduction \cite{Baran1989} and
supply capacity~improvement~\cite{Kening2016}, etc.,
is essentially forming a $\kappa_1$-tree.
That is, a spanning tree is formed if $\kappa_1=1$, while a spanning forest is formed
if $\kappa_1>1$. It is also required that each component has a substation node.

Regarding resilient DS reconfiguration (e.g., MG formation for service restoration),
it is to form a $\kappa_2$-tree
with $\kappa_2 \geq \kappa_1$. In this case, it is required that each component has at most one
substation node, if any.
Actually, the value of $\kappa_2$ should be optimized in resilient reconfiguration.
If it is predefined, the optimality can be impacted.
In all relevant cases, we will~have $\kappa_2 > 1$,
which indicates that a spanning forest is formed.
That is, the topology issues in resilient DS reconfiguration
can be resolved by requiring the network to be a spanning forest. Thus, we can formulate the
radiality constraints as equations \eqref{SF1}-\eqref{SF2}, which are
inspired by the remark as below:

\begin{remark}\label{subgraph}
  An arbitrary subgraph of a spanning tree is~a~spanning forest.
\end{remark}

Specifically, a \emph{subgraph} of a graph consists of a subset of the vertices and edges
in the graph.~For example,~Fig.~\ref{STandSF}(b)~is~a~subgraph of
the graph in Fig.~\ref{STandSF}(a).
According to Definitions~\ref{ST_defi}~and~\ref{SF_defi}, Remark~\ref{subgraph} is naturally true.
Based on that, the proposed radiality constraints are formulated as follows:
\begin{equation}\label{SF1}
  \begin{aligned}
    \boldsymbol{\beta} \in \boldsymbol{\Omega}.
  \end{aligned}
\end{equation}
\begin{equation}\label{SF2}
  \begin{aligned}
    \alpha_{ij} \leq \beta_{ij}, \forall (i,j) \in \boldsymbol{L}.
  \end{aligned}
\end{equation}
In \eqref{SF1}-\eqref{SF2},
$\boldsymbol{L}$ is the set of DS branches;
$\boldsymbol{\alpha}\triangleq\{\alpha_{ij},\forall (i,j)\in \boldsymbol{L}\}$,
where $\alpha_{ij}$ is the connection status of branch $(i,j)$
($1$ if closed, $0$ if open);
$\boldsymbol{\beta} \triangleq \{\beta_{ij},\forall (i,j)\in \boldsymbol{L}\}$,
where $\beta_{ij}$~is~the~\emph{fictitious} connection status of
branch $(i,j)$~($1$~if closed, $0$ if
open).~By ``fictitious'', it indicates that $\boldsymbol{\beta}$ are just auxiliary variables,
which do not actually determine the network topology of the DS.
The topology is still determined by variables $\boldsymbol{\alpha}$.
The symbol~$\boldsymbol{\Omega}$ is the set of all
incidence vectors of spanning tree topologies
that the network can form via reconfiguration (see Section~\ref{tightness_issue}).

Thus, constraint \eqref{SF1} enforces $\boldsymbol{\beta}$ to form a
fictitious spanning tree. Constraint
\eqref{SF2} then restricts the DS to close a subset of the closed branches in the spanning tree
determined by $\boldsymbol{\beta}$.
That is, $\boldsymbol{\alpha}$ form a subgraph of the fictitious spanning tree.
Remark~\ref{subgraph} indicates that the resulting network topology
determined by $\boldsymbol{\alpha}$ is a spanning forest.
Note that constraint \eqref{SF1} is just expressed conceptually here. We will
elaborate on its explicit formulations in Section~\ref{tightness_issue}.
Besides, constraints \eqref{SF1}-\eqref{SF2}
implicitly assume that the DS has only one substation node.
Section~\ref{tightness_issue} will briefly explain the reason for this assumption,
and introduce a simple method that enables the proposed radiality constraints to fully
handle a DS with multiple substation nodes.
{\color{black}For uncontrollable branches without switches or with faulted switches,
constraints can be added to specify their connection status (see Section~\ref{applicationMG}).}

Remark~\ref{subgraph} only specifies that any $\boldsymbol{\alpha}$
satisfying constraints \eqref{SF1}-\eqref{SF2} is topologically feasible for the DS.
The following theorem indicates that any topologically feasible $\boldsymbol{\alpha}$
satisfies constraints \eqref{SF1}-\eqref{SF2}.
As it is somewhat less straightforward, its proof is~also provided.
Remark~\ref{subgraph} and Theorem~\ref{subgraph2} together
certify the validity of the proposed radiality constraints. 

\begin{theorem}\label{subgraph2}
  {\color{black}A} spanning forest subgraph of a connected graph is {\color{black}also
  the} subgraph of {\color{black}a} spanning
  tree~subgraph~of~the~connected graph.
\end{theorem}

\emph{Proof:}
Assume a $\kappa$-tree denoted as $\mathcal{G}_{\kappa}$ to be a subgraph~of a graph $\mathcal{G}$.
As $\mathcal{G}$ is connected, there exists at least one edge in $\mathcal{G}$ that can link the
$\kappa$th component of $\mathcal{G}_{\kappa}$ to another component of $\mathcal{G}_{\kappa}$.
Adding this edge to $\mathcal{G}_{\kappa}$, which does not create cycles, $\mathcal{G}_{\kappa}$
becomes a $(\kappa-1)$-tree denoted as $\mathcal{G}_{\kappa-1}$. Repeating this process,
$\mathcal{G}_{\kappa}$ ends up to be a $1$-tree denoted as $\mathcal{G}_{1}$.
Obviously, the spanning forest $\mathcal{G}_{\kappa}$ is a subgraph of the spanning tree
$\mathcal{G}_{1}$, which is a subgraph of $\mathcal{G}$.
This completes the proof.
\ \ \ \ \ \ \
$\square$

{\color{black}The proposed radiality constraints \eqref{SF1}-\eqref{SF2}
can be interpreted as a two-step method to
regulate DS topology in reconfiguration.
In the first step, constraint \eqref{SF1}
ensures~network~radiality~by~having a spanning tree to be the DS topology's
supergraph, which is determined by
the fictitious connection status of branches. 
(If $\mathcal{G}_{A}$ is a subgraph of $\mathcal{G}_{B}$,
then $\mathcal{G}_{B}$ is said to be a \emph{supergraph} of $\mathcal{G}_{A}$.)
In the second step, constraint \eqref{SF2}
enables more flexible reconfiguration
by allowing the DS
to select a subgraph of the fictitious spanning tree to be the
actual network topology.}

By contrast, common methods of formulating radiality~constraints
in the literature
(e.g., \cite{Lav2012} \cite{Ahm2015})
can be seen as a one-step process
that directly enforces the DS topology to be a~spanning tree or a
spanning forest.
Their models (the single-commodity flow model \cite{Lav2012}, etc.)
can also be used to formulate constraint~\eqref{SF1}
here. Nevertheless, a new formulation will~be presented
in Section~\ref{tightness_issue}.
{\color{black}Above all,
the~proposed~two-step~method~for~radiality constraints formulation
enables~many~more flexibilities~in DS reconfiguration, such as
more adaptive merge or separation~of~sub-grids,
and~more~flexible~allocation~of~power~sources~into sub-grids.
Such benefits will be detailed in
Sections~\ref{applicationMG}~and~\ref{cases_MG}~applying
the proposed radiality constraints to resilient MG formation.}

Due to the critical importance of network reconfiguration~in
DS resilience enhancement, etc., many publications~also~propose to
develop new methods of ensuring radiality to~allow~undiminished
flexibilities and adaptivities of DSs in~optimization. Specifically,
references \cite{SSMA2018}\cite{Arif2018}
have presented two applicable methods, which are compared to
our proposed method as~below:

\subsubsection{Validity}
Remark~\ref{subgraph} and Theorem~\ref{subgraph2}
theoretically prove~the
validity of constraints \eqref{SF1}-\eqref{SF2} here.
References \cite{SSMA2018} \cite{Arif2018} however
lack such analytical proofs for their proposed formulations.

\subsubsection{Tightness}
With constraint \eqref{SF1} properly formulated,
our proposed model attains the tightest formulation
of the \emph{spanning} \emph{forest polytope} of the DS,
i.e., the convex hull of incidence~vectors of
possible spanning forest
topologies (see~Section~\ref{tightness_issue}).
The methods in~\cite{SSMA2018} \cite{Arif2018}
produce
relatively less tight~formulations.

\subsubsection{Compactness}
Constraint~\eqref{SF1} can also be
formulated in less tight but more compact manners (see~Section~\ref{tightness_issue}).
In this regard,
the resulting topology constraints \eqref{SF1}-\eqref{SF2}
are generally more compact than (i.e., with fewer variables and constraints),
and still as tight as or even tighter than those of \cite{SSMA2018} \cite{Arif2018}.

\subsubsection{Application convenience}
The application of our proposed method can be quite straightforward,
i.e., simply adding constraint \eqref{SF2} to the commonly-used
single-commodity flow-based radiality constraints.
The methods in \cite{SSMA2018} \cite{Arif2018} however~involve the
introduction of a virtual source node and~virtual branches,
and the modeling of a virtual DC optimal power~flow subproblem
and its Karush-Kuhn-Tucker conditions etc.,~respectively.

\subsubsection{Applicability}
The proposed radiality constraints \eqref{SF1}-\eqref{SF2} can be
adopted in different DS optimization
problems involving reconfiguration.
Topological and some other related flexibilities thus can be fully enabled in
optimization.
The~methods~in~\cite{SSMA2018} \cite{Arif2018}
have some limitations in this regard.
For example,~the radiality constraints proposed in \cite{SSMA2018}
allow the merge of~sub-grids in optimization,
but do not enable their possible separation.

\section{Tightness and Compactness Issues}\label{tightness_issue}

Constraint \eqref{SF1}
is only expressed conceptually in Section~\ref{proposed_method}.
As aforementioned,
common methods or models~for~representing radiality constraints
in DS reconfiguration-related~publications
can be used for its explicit formulation.
Those methods and models are revisited here.
A new model is~also presented.
The tightness and
compactness issues are specifically
discussed.

DS reconfiguration is essentially a
mixed-integer programming (MIP) problem
generally solved by the branch-and-cut (B\&C) method,
which needs to solve its \emph{linear programming (LP)
relaxations} (i.e., relaxing integer$/$binary constraints).~The
computational complexity of a MIP
depends on its tightness and compactness, etc.
A MIP formulation is \emph{tight} if its feasible region
is similar to that of its LP relaxation,
contributing to a smaller gap between the optimal values (solutions)
of the MIP and its LP relaxation, and less explored nodes in
the B\&C search tree (i.e., fewer iterations and
faster~convergence).~A \emph{compact} MIP formulation has a small number of variables~and constraints,
leading to shorter computation time for~each explored node in
the B\&C search tree.
Tightness and compactness are often conflicting objectives in formulating
a MIP \cite{UCF2013}.

\begin{figure}[t!] 
  \centering
  \includegraphics[width=3.3in]{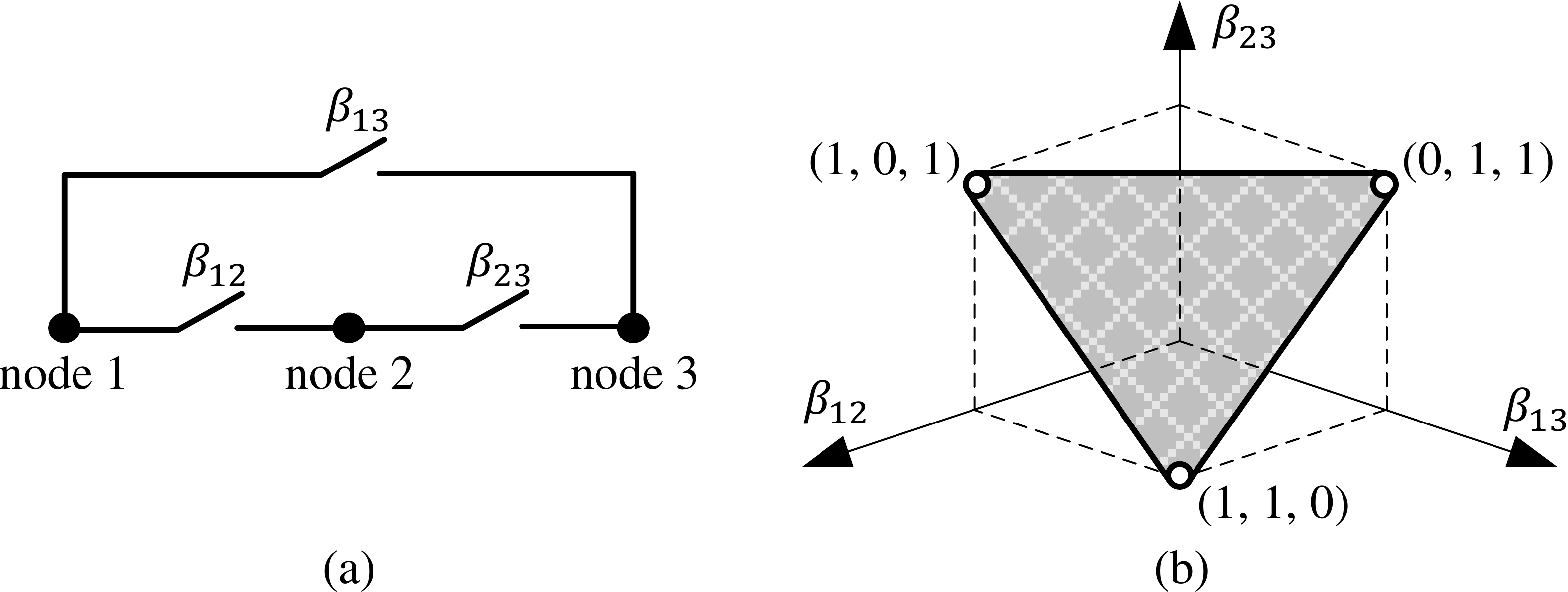}
    \vspace{-8pt}
  \caption{(a) A 3-node sample system; {\color{black}(b) \emph{The set of incidence vectors of~spanning trees}
  (i.e., $\boldsymbol{\Omega}$)
  and \emph{spanning tree polytope} (i.e., convex hull of $\boldsymbol{\Omega}$) illustrated on the sample system.}}
  \label{polytope_eg}
    \vspace{-1pt}
\end{figure}

Several relevant concepts are further introduced.
Here,~\emph{incidence vectors of spanning
trees} are values of $\boldsymbol{\beta}$
defining~fictitious spanning tree topologies of the DS;
and \emph{spanning~tree~polytope} is the convex hull of
such incidence vectors.
For example, the system in Fig.~\ref{polytope_eg}(a)
has $\boldsymbol{\beta}=[\beta_{12},\beta_{13},\beta_{23}]$.
Its incidence vectors of spanning
trees include $[1,1,0]$, $[1,0,1]$ and $[0,1,1]$,
and its spanning tree polytope is the shaded region
with~those vectors as the vertices
in Fig.~\ref{polytope_eg}(b).
\emph{Incidence~vectors~of~spanning forests} and \emph{spanning forest polytope}
are defined likewise.

That is,
constraint \eqref{SF1} actually requires $\boldsymbol{\beta}$
to be an incidence vector of a spanning tree.
A tight and compact~formulation is desired to explicitly represent this requirement.
LP relaxations of the tightest formulations define the spanning tree polytope.

Some methods or models to express radiality constraints~in
DS reconfiguration-related
publications can be used to~formulate the spanning tree constraint~\eqref{SF1}.
They are revisited~here:

\begin{table}[t!]
  \centering
    \vspace{-1pt}
  \caption{Comparing Different Types of Spanning Tree Constraints}
    \vspace{-1pt}
    \setlength{\tabcolsep}{3.3pt}
    \begin{tabular}{|c|c|c|c|}
    \hline
    \tabincell{c}{Type of formu-\\lation/model}
    & \tabincell{c}{Number of \\ variables}
    & \tabincell{c}{Number of constraints}
    & \tabincell{c}{Spanning tree \\polytope ?}  \\
    \hline
    \tabincell{c}{Subtour \\ elimination}
    & $|\boldsymbol{L}|$
    & $\mathcal{O}(2^{|\boldsymbol{N}|})$
    & Yes \\
    \hline
    \tabincell{c}{Directed cutset}
    & $3$$\cdot$$|\boldsymbol{L}|$
    & $\mathcal{O}(2^{|\boldsymbol{N}|-1})$
    & \tabincell{c}{Yes (extended \\formulation)} \\
    \hline
    \tabincell{c}{Single-com-\\modity flow}
    & $3$$\cdot$$|\boldsymbol{L}|$
    & $|\boldsymbol{N}|$$+$$2$$\cdot$$|\boldsymbol{L}|$
    & No \\
    \hline
    \tabincell{c}{Primal \& dual \\graphs-based}
    & $4$$\cdot$$|\boldsymbol{L}|$
    & $2$$\cdot$$|\boldsymbol{L}|$
    & \tabincell{c}{Yes (extended \\formulation)} \\
    \hline
    \tabincell{c}{Directed multi-\\commodity flow}
    & \tabincell{c}{$2$$\cdot$$|\boldsymbol{N}|$$\cdot$$|\boldsymbol{L}|$\\
    $+$$|\boldsymbol{L}|$}
    & \tabincell{c}{$|\boldsymbol{N}|$$^2$$+$$2$$\cdot$$|\boldsymbol{N}|$$\cdot$$|\boldsymbol{L}|$
    \\$-$$|\boldsymbol{N}|$$-$$|\boldsymbol{L}|$$+$$1$}
    & \tabincell{c}{Yes (extended \\formulation)} \\
    \hline
    \end{tabular}%
  \label{STcons}%
    \vspace{-6pt}
\end{table}%

\subsubsection{Loop-eliminating method \cite{Ji1996}}

It enumerates all loops~of the DS, and enforces each one
to be open. However, finding all loops in a graph is NP-hard.
This method is essentially~equivalent to the
\emph{subtour elimination formulation} of spanning tree constraints in
graph theory. As one of the tightest formulations, its LP relaxation
defines the spanning tree polytope and thus has integer extreme points.
It has only $|\boldsymbol{L}|$ variables but also~an exponential number of
constraints, which limits its application.

\subsubsection{Path-based model \cite{Ram2005}}
It enumerates~all~paths~to~the~substation for each node,
and activates only one of them.~Still,~finding all paths between two graph nodes
is
NP-hard.~This~model also has a counterpart in graph theory,
i.e., the \emph{directed~cutset formulation} of spanning tree constraints with
$3$$\cdot$$|\boldsymbol{L}|$~variables and~an exponential number of
constraints. It is one of the tightest formulations
with integer vertices of its LP relaxation,~too.

\subsubsection{Single-commodity flow-based model \cite{Lav2012}}

It closes $|\boldsymbol{N}|$$-$$1$ branches ($\boldsymbol{N}$: set of DS nodes),
and ensures connectivity~by~imposing unit fictitious flow from the substation to
each node.
With $3$$\cdot$$|\boldsymbol{L}|$ (or reduced to
$2$$\cdot$$|\boldsymbol{L}|$) variables and a linear number of constraints,
it is one of the most compact formulations.
Due to its simplicity,
it is also the~most~commonly~used~DS~radiality model.
Still, it is less tight. The projection of its~LP~relaxation
into the $\boldsymbol{\beta}$-space defines a
region larger than the~spanning tree polytope, and therefore has fractional
extreme points.

\subsubsection{Parent-child relation-based method \cite{Jab2012}}

It~instructs~every node except the substation to have one parent.
However, this method can produce a disconnected graph with loops \cite{Ahm2015}.

\subsubsection{Primal and dual graphs-based model \cite{Ahm2015}}

It~froces~each primal (dual) node to be connected to another primal (dual)
node, and forbids the primal and dual spanning trees to~be~intersected.
With $4$$\cdot$$|\boldsymbol{L}|$ variables
and a linear number of constraints,
it is one of the most compact formulations.
Originally~proposed in \cite{Williams2002},
it actually has a totally unimodular constraint matrix
to define a polyhedron with integer vertices.
Thus, it~is~also~one~of the tightest formulations with its
LP relaxation~defining~the spanning tree polytope
in the $\boldsymbol{\beta}$-space.
Still, it is only applicable to planar graphs,
and needs to construct dual graphs.
It~also~has a flaw to be avoided by
appropriate selection~of~roots~\cite{Validi2019}.~Generally,
its applicability is limited and its use is not convenient.

\begin{table}[t!]
  \centering
  \caption{Notations of the Spanning Tree Constraints \eqref{dflo_1}-\eqref{dflo_7}}
    \vspace{-6pt}
    \begin{tabular}{|c|l|}
    \hline
    \multicolumn{2}{|c|}{\emph{Parameters}} \\
    \hline
    $\boldsymbol{N}/\boldsymbol{L}$ & Set of all DS nodes/branches. \\
    \hline
    $i_r$ & Index for the substation node. \\
    \hline
    \hline
    \multicolumn{2}{|c|}{\emph{Variables}} \\
    \hline
    $f_{ij}^k$    &  Flow of commodity $k$ from node $i$ to node $j$.  \\
    \hline
    $\lambda_{ij}$ &  \tabincell{l}{Binary, $1$ if arc $(i,j)$ is included
    in the \\directed spanning tree, $0$ otherwise.} \\
    \hline
    $\beta_{ij}$ & \tabincell{l}{Binary, $1$ if the fictitious connection
    status of \\ branch $(i,j)$ is closed, $0$ if open.}
    \\
    \hline
    \end{tabular}
  \label{notations_zero}%
    \vspace{-11pt}
\end{table}%

Using notations in Table~\ref{notations_zero},
we further introduce a new~flow-based formulation of the spanning tree constraint
\eqref{SF1}
as below:
\begin{equation}\label{dflo_1}
  \begin{aligned}
    \sum_{(j,i_r)\in\boldsymbol{L}} f_{ji_r}^k
    - \sum_{(i_r,j)\in\boldsymbol{L}} f_{i_rj}^k = -1, \forall k \in
    \boldsymbol{N} \setminus i_r.
  \end{aligned}
\end{equation}
\begin{equation}\label{dflo_2}
  \begin{aligned}
    \sum_{(j,i)\in\boldsymbol{L}} f_{ji}^k
    - \sum_{(i,j)\in\boldsymbol{L}} f_{ij}^k = 0,
    \forall k \in \boldsymbol{N} \setminus i_r,
    \forall i \in \boldsymbol{N} \setminus \{i_r,k\}.
  \end{aligned}
\end{equation}
\begin{equation}\label{dflo_3}
  \begin{aligned}
    \sum_{(j,k)\in\boldsymbol{L}} f_{jk}^k
    - \sum_{(k,j)\in\boldsymbol{L}} f_{kj}^k = 1,
    \forall k \in \boldsymbol{N} \setminus i_r.
  \end{aligned}
\end{equation}
\begin{equation}\label{dflo_4}
  \begin{aligned}
    0 \le f_{ij}^k \le \lambda_{ij}, 0 \le f_{ji}^k \le \lambda_{ji},
    \forall k \in \boldsymbol{N} \setminus i_r,
    \forall (i,j) \in \boldsymbol{L}.
  \end{aligned}
\end{equation}
\begin{equation}\label{dflo_5}
  \begin{aligned}
    \sum_{(i,j)\in\boldsymbol{L}} (\lambda_{ij}+\lambda_{ji})
    = |\boldsymbol{N}| - 1.
  \end{aligned}
\end{equation}
\begin{equation}\label{dflo_6}
  \begin{aligned}
    \lambda_{ij}+\lambda_{ji} = \beta_{ij},
    \forall (i,j) \in \boldsymbol{L}.
  \end{aligned}
\end{equation}
\begin{equation}\label{dflo_7}
  \begin{aligned}
    \lambda_{ij},\lambda_{ji} \in \{0,1\}, \forall (i,j) \in \boldsymbol{L}.
  \end{aligned}
\end{equation}
The above formulation is called as the
\emph{directed multicommodity flow-based model}
of the spanning tree constraints.
It~defines~a \emph{fictitious} commodity for
each node $k\ne i_r$, and enforces 1 unit of commodity $k$ delivered
from the substation node $i_r$ to node $k$. Constraint \eqref{dflo_4} implies
that each commodity can flow on an arc only if the arc is included in the
directed spanning tree defined by
variables $\lambda_{ij}$. Other equations are self-explanatory.

The LP relaxation of the above formulation \eqref{dflo_1}-\eqref{dflo_7}
defines the spanning tree polytope
(i.e., $\mathit{conv}(\boldsymbol{\Omega})$)
in the $\boldsymbol{\beta}$-space
by a polynomial number of variables and constraints.
As compared in Table~\ref{STcons},
it is generally the most compact formulation
among the ones that are the tightest and applicable for
both planar and non-planar graphs.
Now the models listed in Table~\ref{STcons}
actually cover a wide spectrum of tightness and compactness levels
for the formulation of spanning tree constraints.
Researchers and practitioners may choose an appropriate model based on their
preferences and needs, etc.

Next,
the tightness and compactness of
constraints~\eqref{SF1}-\eqref{SF2},
rather than solely constraint~\eqref{SF1},
are further discussed briefly.

\begin{proposition}\label{SFP}
  If the LP relaxation of the explicit formulation of constraint \eqref{SF1}
  defines the
  spanning tree polytope~in~the~$\boldsymbol{\beta}$-space,
  the LP relaxation of the proposed radiality constraints \eqref{SF1}-\eqref{SF2}
  defines the
  spanning forest polytope in the~$\boldsymbol{\alpha}$-space.
\end{proposition}

\emph{Proof:}
Assume that
the LP relaxation of constraints \eqref{SF1}-\eqref{SF2} has
a fractional vertex
$(\boldsymbol{\alpha^*},\boldsymbol{\beta^*})$
in the $(\boldsymbol{\alpha},\boldsymbol{\beta})$-space.~Further assume that
$\boldsymbol{\beta^*}$ has fractional entries.
Thus, $\boldsymbol{\beta^*}$
can be~represented by an integer vertex
$\boldsymbol{\beta^m}$
and another point $\boldsymbol{\beta^n}$ of~the spanning tree polytope:
$\boldsymbol{\beta^*}=\xi\cdot\boldsymbol{\beta^m}+(1-\xi)\cdot\boldsymbol{\beta^n}$
with~$0<\xi<1$.
For $(i,j)$ 
with $\beta_{ij}^{m}=0$ or
$\beta_{ij}^{m}=1$,
we have
$\beta_{ij}^{*}=(1-\xi)\cdot\beta_{ij}^{n}$
and $\beta_{ij}^{*}=\xi+(1-\xi)\cdot\beta_{ij}^{n}$, respectively.
As $0\le\beta_{ij}^{n}\le1$ and $0\le\alpha_{ij}^{*}$ $\le\beta_{ij}^{*}\le1$,
we have $0\le\frac{\alpha_{ij}^{*}}{1-\xi}\le\beta_{ij}^{n}$
for $(i,j)$ with $\beta_{ij}^{m}=0$,
and have $0\le\frac{\alpha_{ij}^{*}-\xi\cdot\chi_{ij}}{1-\xi}\le\beta_{ij}^{n}$
with
$\chi_{ij}\in[\frac{\alpha_{ij}^{*}-\beta_{ij}^{*}+\xi}{\xi},\frac{\alpha_{ij}^{*}}{\xi}]\cap[0,1]$
for
$(i,j)$ with $\beta_{ij}^{m}=1$.
Then, we can construct
$\boldsymbol{\alpha^m}$ and~$\boldsymbol{\alpha^n}$~to have
$\alpha_{ij}^{m}=0$ and $\alpha_{ij}^{n}=\frac{\alpha_{ij}^{*}}{1-\xi}$
for $(i,j)$ with $\beta_{ij}^{m}=0$,
and have $\alpha_{ij}^{m}=\chi_{ij}$ and $\alpha_{ij}^{n}=\frac{\alpha_{ij}^{*}-\xi\cdot\chi_{ij}}{1-\xi}$
for $(i,j)$ with $\beta_{ij}^{m}=1$,
so~that $(\boldsymbol{\alpha^*},\boldsymbol{\beta^*})=\xi\cdot(\boldsymbol{\alpha^m},\boldsymbol{\beta^m})+(1-\xi)\cdot(\boldsymbol{\alpha^n},\boldsymbol{\beta^n})$.
This~contradicts the assumption of $(\boldsymbol{\alpha^*},\boldsymbol{\beta^*})$
being a vertex.
The case with $\boldsymbol{\beta^*}$ as a vertex of $\mathit{conv}(\boldsymbol{\Omega})$
can be analyzed similarly.~Thus,~vertices of the projection of
constraints \eqref{SF1}-\eqref{SF2}'s LP relaxation~into $\boldsymbol{\alpha}$-space are~0-1
incidence vectors of spanning forests of the~DS.~$\square$

Proposition~\ref{SFP} implies that
the tightness and compactness~features
of constraints \eqref{SF1}-\eqref{SF2}
essentially~follow the~explicit~formulation
of constraint~\eqref{SF1}.
Thus, the proposed radiality~constraints \eqref{SF1}-\eqref{SF2} also
cover a wide spectrum of
tightness and compactness levels for the formulation of spanning
forest polytope.~An~appropriate model can be selected
based~on~one's~preferences,~etc.

Specifically,
to exploit the aforementioned properties
of the formulations,
the DS is assumed to have only one substation node.
For a DS with multiple substation nodes,
one can~merge them into one node in modeling constraints
\eqref{SF1}-\eqref{SF2}, but still treat them as separate nodes in
modeling DS operational constraints.

\section{Application to Resilient MG Formation}\label{applicationMG}

The proposed radiality model \eqref{SF1}-\eqref{SF2}
can be applied in~different optimization problems involving
DSs and/or MGs with flexible boundaries \cite{Arefifar2015, Kim2018, Manshadi2018}, etc.
Here, we apply it~to~regulate~the DS topology in resilient MG formation
to~verify~its~advantages.

\begin{table}[t!]
  \centering
  \caption{Notations of the MG Formation Model \eqref{obj}-\eqref{PF13}}
  \vspace{-3pt}
    \setlength{\tabcolsep}{1.38pt}
    \begin{tabular}{|c|l|}
    \hline
    \multicolumn{2}{|c|}{\emph{Parameters}} \\
    \hline
    $\boldsymbol{N_{r}}/\boldsymbol{N_{g}}$
    &
    Set of substation nodes/DG nodes.
    \\
    \hline
    $\boldsymbol{N_{o}}/\boldsymbol{N_{c}}$ &
    Set of nodes with faulted open/closed load switches.\\
    \hline
    $\boldsymbol{L_{o}}/\boldsymbol{L_{c}}$ &
    Set of faulted open/closed branches.
    \\
    \hline
    $p_{i}^{c}/q_{i}^{c}$     &    Real/reactive power demand of the load at node $i$.  \\
    \hline
    $\overline{p}_{i}^{g}/\overline{q}_{i}^{g}$ & \tabincell{l}{Real/reactive power capacity of the power source at node $i$.} \\
    \hline
    $\overline{v}_{i}/\underline{v}_{i}$    &     \tabincell{l}{Maximum/minimum squared voltage magnitude of node $i$.} \\
    \hline
    $r_{ij}/x_{ij}/\overline{S}_{ij}$     &   Resistance/reactance/apparent power capacity of branch $(i,j)$. \\
    \hline
    $\omega_i$ & Priority weight of the load at node $i$.\\
    \hline
    $\vartheta_i$ &  Number of branches starting or ending with node $i$.\\
    \hline
    $M$ & A large enough positive number. \\
    \hline
    \hline
    \multicolumn{2}{|c|}{\emph{Variables}} \\
    \hline
    $\delta_{i}$  & \tabincell{l}{Binary, $1$ if the load at node $i$ is picked up, $0$ otherwise.}  \\
    \hline
    $\varepsilon_i$ & Binary, $1$ if node $i$ is energized, $0$ otherwise. \\
    \hline
    $\alpha_{ij}$ & Binary, $1$ is branch $(i,j)$ is closed, $0$ if open. \\
    \hline
    $p_{i}^{g}/q_{i}^{g}$    & \tabincell{l}{Real/reactive power output of the power source at node $i$.} \\
    \hline
    $v_{i}$     &   Squared voltage magnitude of node $i$.\\
    \hline
    $P_{ij}/Q_{ij}$   &    Real/reactive power flow on branch $(i,j)$.\\
    \hline
    \end{tabular}%
  \label{notations}%
    \vspace{-6pt}
\end{table}%

MG formation has recently been
extensively studied~in~\cite{Bie2017}
\cite{Chen2016} \cite{Ding2017} \cite{Ding2017_2} \cite{Sedzro2018}, etc.
It is to reconfigure the DS~into~multiple MGs energized
by DGs, so as to restore critical loads,~etc.~The
literature currently has two major types of optimization models
for this problem \cite{Chen2016} \cite{Ding2017}. Using the proposed
radiality model, a new formulation for resilient MG formation
is constructed:
\begin{equation}\label{obj}
  \begin{aligned}
    \max \sum_{i \in \boldsymbol{N}} \delta_i \cdot \omega_i \cdot p_{i}^{c}
  \end{aligned}
\end{equation}
\begin{equation*}\label{MG_SF}
  \begin{aligned}
    s.t.\ \ \ \ \ \ \ \ \ \eqref{SF1}-\eqref{SF2}. \ \ \ \ \ \ \ \ \ \ \ \ \
  \end{aligned}
\end{equation*}
\begin{equation}\label{PF1}
  \begin{aligned}
    p_{i}^{g} - \delta_{i} \cdot p_{i}^{c}
    + \sum_{(j,i) \in \boldsymbol{L}} P_{ji}
    - \sum_{(i,j) \in \boldsymbol{L}} P_{ij} = 0,
    \forall i \in \boldsymbol{N}.
  \end{aligned}
\end{equation}
\begin{equation}\label{PF2}
  \begin{aligned}
    q_{i}^{g} - \delta_{i} \cdot q_{i}^{c}
    + \sum_{(j,i) \in \boldsymbol{L}} Q_{ji}
    - \sum_{(i,j) \in \boldsymbol{L}} Q_{ij} =0,
    \forall i \in \boldsymbol{N}.
  \end{aligned}
\end{equation}
\begin{equation}\label{PF3}
  \begin{aligned}
    p_{i}^{g} = q_{i}^{g} = 0, \forall i \in \boldsymbol{N} \setminus
    \{\boldsymbol{N_{r}}, \boldsymbol{N_{g}}\}.
  \end{aligned}
\end{equation}
\begin{equation}\label{PF3.5}
  \begin{aligned}
    0 \le p_{i}^{g} \le \overline{p}_{i}^{g},0 \le q_{i}^{g} \le \overline{q}_{i}^{g},
    \forall i \in \{\boldsymbol{N_{r}}, \boldsymbol{N_{g}}\}.
  \end{aligned}
\end{equation}
\begin{equation}\label{PF4}
  \begin{aligned}
    v_i - v_j \ge 2 \cdot (P_{ij} \cdot r_{ij} + Q_{ij} \cdot x_{ij})
    + (\alpha_{ij} - 1) \cdot M,&
    \\
    v_i - v_j \le  2 \cdot (P_{ij} \cdot r_{ij}
    + Q_{ij} \cdot x_{ij})
    + (1 - \alpha_{ij}) \cdot M,&
    \\
    \forall (i,j) \in \boldsymbol{L}.&
  \end{aligned}
\end{equation}
\begin{equation}\label{PF6}
  \begin{aligned}
    \underline{v}_i \le {v}_i \le \overline{v}_i, \forall i \in \boldsymbol{N}.
  \end{aligned}
\end{equation}
\begin{equation}\label{PF7}
  \begin{aligned}
    P_{ij}^{2} + Q_{ij}^{2} \le \alpha_{ij} \cdot \overline{S}_{ij}^2, \forall (i,j) \in \boldsymbol{L}.
  \end{aligned}
\end{equation}
\begin{equation}\label{PF9}
  \begin{aligned}
    \alpha_{ij} = 0, \forall (i,j) \in \boldsymbol{L_{o}};
      \alpha_{ij} = 1, \forall (i,j) \in \boldsymbol{L_{c}}.
  \end{aligned}
\end{equation}
\begin{equation}\label{PF11}
  \begin{aligned}
    \delta_{i} = 0, \forall i \in \boldsymbol{N_{o}};
    \delta_{i} \ge \varepsilon_i, \forall i \in \boldsymbol{N_{c}}.
  \end{aligned}
\end{equation}
\begin{equation}\label{PF12}
  \begin{aligned}
    \varepsilon_i = 1, \forall i \in \{\boldsymbol{N_{r}}, \boldsymbol{N_{g}}\}.
  \end{aligned}
\end{equation}
\begin{equation}\label{PF13}
  \begin{aligned}
   (\sum_{(i,j)\in \boldsymbol{L}} \varepsilon_j \cdot \alpha_{ij}
    + \sum_{(j,i)\in \boldsymbol{L}} \varepsilon_j \cdot \alpha_{ji})/
    \vartheta_i
    \le \varepsilon_i
    \ \ \ \ \ \ \ \ \ \ \ \ \ \ \ \ \ \ \ \ \ \ \
\\    \le
    \sum_{(i,j)\in \boldsymbol{L}} \varepsilon_j \cdot \alpha_{ij}
    + \sum_{(j,i)\in \boldsymbol{L}} \varepsilon_j \cdot \alpha_{ji},
    \forall i \in
    \boldsymbol{N} \setminus
    \{\boldsymbol{N_{r}}, \boldsymbol{N_{g}}\}.
  \end{aligned}
\end{equation}

Notations are listed in Table~\ref{notations}.
The objective~function~\eqref{obj}
maximizes the weighted sum of restored loads.~Constraints~\eqref{SF1}-\eqref{SF2}
are used to ensure radiality and enable topological~flexibilities,
etc.
{\color{black}Constraint \eqref{SF1} is formulated
using the methods in Section~\ref{tightness_issue}.}
Equations \eqref{PF1}-\eqref{PF2} enforce real and~reactive~power~balance, respectively.
Equation \eqref{PF3} imposes zero power output
for nodes without power sources.
Constraint~\eqref{PF3.5}~indicates~real and reactive power capacities of substations or DGs.
Constraint \eqref{PF4} represents the
DistFlow~model~with~the~much~smaller~quadratic terms ignored \cite{Baran1989} \cite{Taylor2012}.
It is relaxed for open branches.
Constraint \eqref{PF6} expresses the voltage magnitude limits.
Constraint \eqref{PF7} is
convex though non-linear, and can be linearized by the technique in \cite{ChenXin2016}, etc.
It limits the apparent power on a branch by its capacity.
Equation \eqref{PF9} restricts
the connection status
of faulted open or faulted
closed branches.
Equation \eqref{PF11}
prohibits picking up the loads with faulted open switches,
and enforces picking up the loads with faulted closed switches
as long as their located nodes are energized.
Equation \eqref{PF12} specifies
that substation and DG nodes are energized.
Constraint \eqref{PF13} derives
the energization status of other nodes by examining
if they are connected to an energized node.
The non-linear and non-convex terms, i.e., $\varepsilon_j\cdot\alpha_{ij}$
and $\varepsilon_j\cdot\alpha_{ji}$,
can be~equivalently~linearized and convexified by the
McCormick envelopes \cite{McCormick1976}.

Note that the areas with surviving access to the main~grid power via substation nodes
are also considered.~Operational constraints may prohibit those areas
to be fully restored by the substations, and forming MGs with DGs
can help achieve better restoration of them \cite{Lei2016}. For statement simplicity,
a sub-grid powered
by a substation is also counted~as~a~MG~here.

\begin{table}[t!]
  \centering
  \caption{Comparing Different MG Formation Models}
    \vspace{-2pt}
    \setlength{\tabcolsep}{4.25pt}
    \begin{tabular}{|c|c|c|c|}
    \hline
    --------- & Model in~\cite{Chen2016} & Model in~\cite{Ding2017} & Proposed model\\
    \hline
    \tabincell{c}{Applicable \\systems} & Radial &
    \tabincell{c}{Radial \\ or meshed} &
    \tabincell{c}{Radial \\ or meshed} \\
    \hline
    \tabincell{c}{Radiality \\constraints} & ------
    & \tabincell{c}{Single-commo- \\dity flow-based}
    & \tabincell{c}{Proposed radi- \\ality constraints} \\
    \hline
    \tabincell{c}{Allocation of \\DGs into MGs}
    & \tabincell{c}{One DG \\in each MG}
    & \tabincell{c}{One DG \\in each MG}
    & Flexible \\
    \hline
    Number of MGs & Fixed & Fixed & Flexible \\
    \hline
    \tabincell{c}{Unenergized islands}
    & Not allowed
    & Not allowed
    & Allowed \\
    \hline
    \tabincell{c}{Loads w/ faulted \\ closed switches} &
    \tabincell{c}{Forced to \\ pick up} &
    \tabincell{c}{Forced to \\ pick up} &
    Flexible \\
    \hline
    \end{tabular}%
  \label{comparing}%
    \vspace{-6pt}
\end{table}%

Table~\ref{comparing} compares the proposed model with the two existing MG
formation models in the literature. The summaries~in~the~table are self-explanatory.
Brief explanations are given as~below:

1) The model in \cite{Chen2016} is designed for radial systems.~It~is~extended in \cite{Sedzro2018} to deal with meshed systems.
Still,~the~extended model does not eliminate loops.
Thus, it may be used for~the bulk power system,
but does not apply to meshed DSs,~which have to open some branches to
be~operated~in radial topologies.
Both model in \cite{Ding2017} and our model can handle~meshed~DSs,
as radiality constraints are included to avoid loops when reconfiguring the network.
Single-commodity flow-based radiality constraints
are used in\cite{Ding2017}, while
our model adopts the proposed radiality constraints
to~fully~enable~topological~flexibilities, etc.

2) Both models in~\cite{Chen2016} and~\cite{Ding2017} allocate one DG
to each MG. Therefore, the number of MGs is actually
fixed and equal to the number of DGs and/or other power sources.
As for our proposed model, the number of DGs designated to different MGs,
and the resulting number of MGs, are flexible.
Thus, it can also be used for dynamic MG formation,
which involves adaptive merge and/or separation of MGs when damaged parts of
the DS
are sequentially repaired, etc.
Such flexibilities introduce many benefits. For example,
larger MGs can~be~formed~to~better match DGs
with different-sized loads, so as to enhance capacity utilization rates
of DGs and the restoration of critical loads.

3) Both models in \cite{Chen2016} \cite{Ding2017}
energize all nodes included~in $\boldsymbol{N}$.
Thus, their models
cannot consider the nodes in load islands without power sources
and isolated by faulted open~branches.
Discrete/non-dispatchable loads
at nodes $\boldsymbol{N_c}$ 
are also forced to be picked up in their models.
Our model can consider~such load islands,
and can be easily extended to optimize the~allocation of mobile power sources in
such islands and their~merge~with other islands after faulted open branches
are repaired.
Our model can also intentionally form unenergized
islands in~optimization, enabling more flexible pick-up of the loads at nodes $\boldsymbol{N_c}$.
This flexibility
can be critical, as some DSs may have a large portion of
loads not equipped with switches, causing $\boldsymbol{N_c}$ to be a large set.
In \cite{Bie2017}, with additional binary variables, the model in \cite{Ding2017} is modified to
permit more flexible allocation of DGs into MGs.
However, it does not fully enable topological and some related flexibilities.
For example, it still requires energizing all nodes.

In general, as indicated in Table~\ref{comparing},
the proposed~MG~formation model is more adaptive and allows
more flexibilities.

\begin{figure}[t!] 
  \centering
  \includegraphics[width=3.3in]{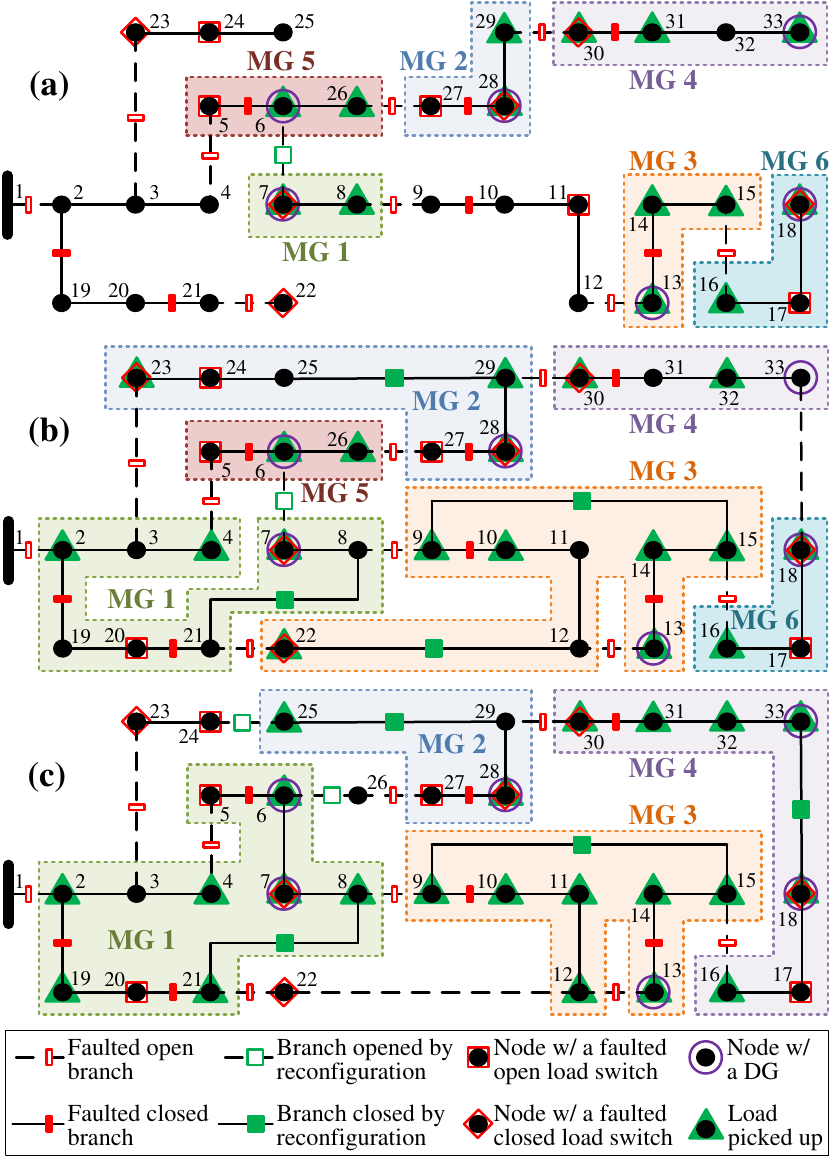}
  \vspace{-4pt}
  \caption{An illustrative case of MG formation on IEEE 33-node test system:
  (a)$/$(b)~Use the model in \cite{Chen2016}$/$\hspace{1sp}\cite{Ding2017};
  (c) Use our proposed~model.}
  \label{33node_base_case}
    \vspace{-6pt}
\end{figure}

\section{Case Studies}\label{cases_MG}

In this section, the proposed resilient MG formation model is
demonstrated on two test systems.
We use a computer with an Intel i5-4278U processor and 8GB of memory.
Involved MIP problems are solved by Gurobi 7.5.2
with the~default~settings.

\subsection{IEEE 33-Node Test System \cite{Baran1989}}

Fig.~\ref{33node_base_case} depicts an illustrative case
based on a~scenario~with~27

\noindent
faults in total.
The model in \cite{Chen2016}
does not apply to~meshed~DSs.
Thus, as indicated in Fig.~\ref{33node_base_case}(a),
it does not consider the normally
open~branches.~It~forms~6~MGs~energized~by~the~6~DGs,~respectively.
As the model in \cite{Ding2017} can handle
meshed~DSs,~it~gets~better results via
further reconfiguration involving those normally open switches, etc.
Fig.~\ref{33node_base_case}(b) shows that it
restores~more~loads~by forming
6 larger MGs, which is essentially a 6-tree.
By~contrast, as in Fig.~\ref{33node_base_case}(c),
our~proposed~model~forms~a~7-tree, i.e., 4 MGs and 3 load islands.
Specifically,~MG~4 and MG~6~in Fig.~\ref{33node_base_case}(b) are
merged into a single MG in Fig.~\ref{33node_base_case}(c),
so~that~the loads
at nodes 31 and 33 can also be restored.
MG 2 in Fig.~\ref{33node_base_case}(b) is~separated into
a MG and a load island in Fig.~\ref{33node_base_case}(c),
so that node~23~is~not energized and its load is not forced to be picked up.~For~space limit, we do not
detail~on~other~differences.~The~models in \cite{Chen2016} \cite{Ding2017}
and our proposed model
have DG capacity~utilization~rates of 55.6\%, 67.4\% and 95.4\%, respectively.
They
restore~1500kW, 1820kW and 2575kW loads, respectively.
Our model achieves more coordinated matching among the different-sized DGs and loads,
and thus attains better service restoration.
Generally, as the proposed radiality constraints can
fully enable topological and some other related flexibilities of the DS,
our MG formation model has extended feasibility and enhanced optimality.

\begin{table}[t!]
  \centering
  \caption{Computation Time and the Number of Infeasible Cases (IEEE 33-Node Test System)}
    \vspace{-6pt}
    \setlength{\tabcolsep}{5.8pt}
    \begin{tabular}{|c|c|c|c|}
    \hline
    --------- & Model in \cite{Chen2016} & Model in \cite{Ding2017} & Proposed model \\
    \hline
    Avg. comput. time & 0.56 s & 0.75 s & 0.49 s \\
    \hline
    Infeasible cases &
    212/10000 &
    272/10000  &
    2/10000 \\
    \hline
    \end{tabular}%
  \label{infeasible_number}%
    \vspace{-5pt}
\end{table}%

\begin{table}[t!]
  \centering
  \caption{Summary Statistics for the Restored Loads (IEEE 33-Node Test System)}
    \vspace{-6pt}
  \setlength{\tabcolsep}{2pt}
    \begin{tabular}{|c|c|c|c|c|c|c|}
    \hline
    \multirow{2}[0]{*}{------} & \multicolumn{5}{c|}{Restored loads (kW)}      & Box plots \\
    \cline{2-7}
     & avg.  & std.  & max.  & med. & min.  & \multirow{4}[0]{*}{\raisebox{-0.1in}{\includegraphics[width=0.95in]{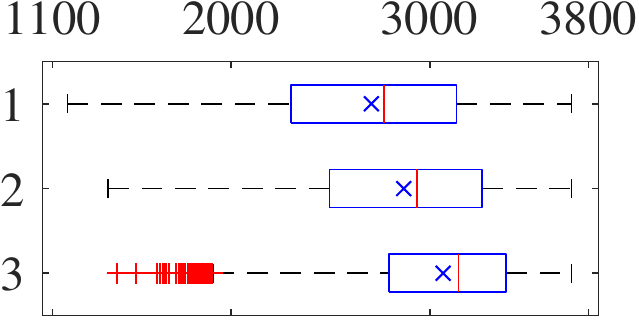}}} \\
     \cline{1-6}
    \tabincell{c}{1: Model in \cite{Chen2016}}     &  2686  &  531  &    3715   &   2750    &   1175    &  \\
    \cline{1-6}
    \tabincell{c}{2: Model in \cite{Ding2017}}     &  2849 &  495  &   3715    &    2915   &    1380   &  \\
    \cline{1-6}
    \tabincell{c}{3: Proposed model}     &   3088    &   393   & 3715 &   3165   &   1425    &  \\
    \hline
    \end{tabular}%
  \label{box_obj}%
    \vspace{-5pt}
\end{table}%

\begin{table}[t!]
  \centering
  \caption{Summary Statistics for the DG Capacity Utilization Rate (IEEE 33-Node Test System)}
    \vspace{-6pt}
  \setlength{\tabcolsep}{0.99pt}
    \begin{tabular}{|c|c|c|c|c|c|c|}
    \hline
    \multirow{2}[0]{*}{------} & \multicolumn{5}{c|}{DG capacity utilization rate}      & Box plots \\
    \cline{2-7}
     & avg.  & std.  & max.  & med. & min.  & \multirow{4}[0]{*}{\raisebox{-0.1in}{\includegraphics[width=1.02in]{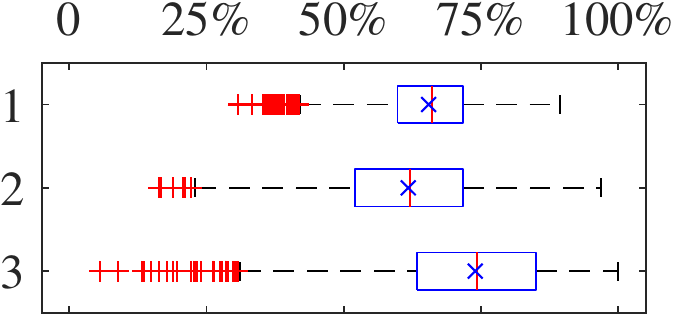}}} \\
     \cline{1-6}
    \tabincell{c}{1: Model in \cite{Chen2016}}     &  64.4\%  &  9.8\%  &    88.4\%  &   65.1\%    &   29.7\%    &  \\
    \cline{1-6}
    \tabincell{c}{2: Model in \cite{Ding2017}}     &  60.7\% &  15.1\%  &   94.8\%    &    61.0\%   &    15.3\%   &  \\
    \cline{1-6}
    \tabincell{c}{3: Proposed model}     &   74.1\%    &   15.7\%   & 100\% &   74.4\%   &   6.6\%    &  \\
    \hline
    \end{tabular}%
  \label{box_DG}%
    \vspace{-7pt}
\end{table}%

\begin{figure}[t!]
\hspace*{\fill}%
  \begin{minipage}[t!]{4.4cm}
    \centering
   \vspace{0pt}
    \includegraphics[height=1.45in]{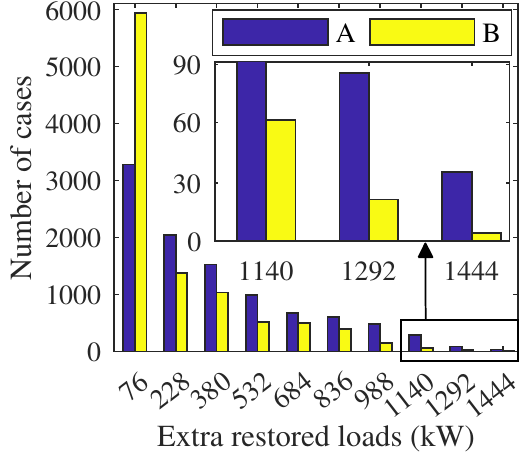}
    \label{extra_loads}
  \end{minipage}%
\hfill
  \begin{minipage}[t!]{4.4cm}
    \centering
    \vspace{0pt}
    \includegraphics[height=1.45in]{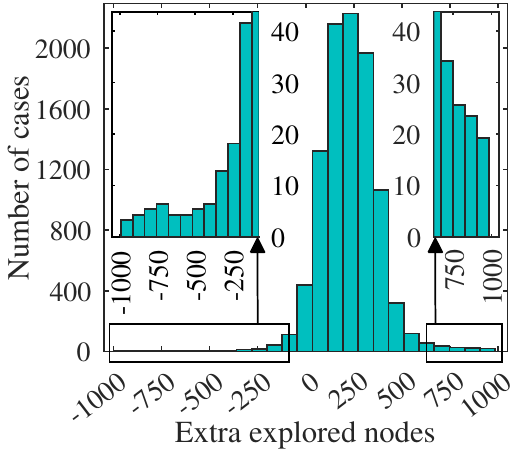}
    \label{extra_nodes}
  \end{minipage}
  \hspace*{\fill}
{\hfill}
  \vspace{-9pt}
  \caption{Left: Histograms of the extra restored loads of our proposed~MG~formation
  model, compared with the model in \cite{Chen2016} (A),~and~compared~with~the~model in \cite{Ding2017} (B).
  Right: Histogram of the extra~explored nodes in the B\&C search tree of
  the proposed MG formation model~using~the~single-commodity flow-based method to formulate
  constraint \eqref{SF1}, compared~with~that~using~the~directed~multicommodity
  flow-based method. (IEEE 33-node~test~system.)}
  \label{histograms}
  \vspace{-2pt}
\end{figure}

To establish superiority of the proposed
radiality constraints and MG formation model,
we further run 10000~cases~based~on randomly generated scenarios of DS faults.
Table~\ref{infeasible_number} shows
that our model has far less infeasible cases.
In many scenarios,~it finds a feasible operating point of
the DS, while the~models~in \cite{Chen2016} \cite{Ding2017}
return infeasibility.
Such results verify that our model
has extended feasibility as the proposed radiality constraints
fully enable topological and some related flexibilities.
Besides,
while its search space is the largest,
its average computation time is the shortest.
Thus, its enlarged feasible set
may possess more computationally tractable characteritics.
As in Table~\ref{box_obj},
our model has the highest average/median/minimum~restored loads
and the smallest standard deviation.
The same maximum value corresponds to the cases with
all loads restored.~On~average, the proposed model restores 15.0\% and 8.4\% more~loads
than
models in \cite{Chen2016}
and \cite{Ding2017},
respectively.~Actually,~our~model performs equally well or better in all cases.
The histograms in Fig.~\ref{histograms}(Left)
depict the outperformance.
Such results validate the enhanced optimality of the proposed model.
Here, one of the
major reasons for its superiority is
more coordinated matching
among~the~DGs~and~loads.~Table~\ref{box_DG}~indicates~that~the~average
DG capacity utilization rate of
our model is much higher than those of the models in \cite{Chen2016} \cite{Ding2017}.
The smaller minimum~and larger standard deviation
are due to the cases with the
substation
contributing much power injection for service~restoration.

The single-commodity flow-based model is used to formulate
constraint \eqref{SF1} in all previous cases. Here, for comparison,~the~directed
multicommodity flow-based model is~used~instead.~Consequently,
the proposed radiality constraints and MG formation
model become tighter, though less compact. The revised model is run
on the same 10000 scenarios.
As indicated in Table \ref{node_time},
$N_m$ is much smaller than
$N_s$, both on average and in 95.8\% of the cases.
The histogram in Fig. \ref{histograms}(Right) also
details the extra explored nodes in the B\&C search tree of the
less tight model.
Although $T_s$ is slightly shorter than $T_m$ on average,
$T_m$ is much shorter than $T_s$ in 10.9\% cases.
Specifically,
in those cases,
$T_s$ and $T_m$ are 1.92 s and 0.69 s on average, respectively;
$N_s$ and $N_m$ are 815 and 186 on average, respectively.
In~general,~the
computation time of the revised tighter model
is more consistent. It also reduces the computation time in
many cases that require the less tight model to explore
much more nodes in the B\&C search tree.
That is, the proposed radiality constraints
can satisfy
the need for tighter formulations of
DS reconfiguration-related optimization problems.

\begin{table}[t!]
  \centering
  \vspace{-1pt}
  \caption{Computation Time and the Number of Explored Nodes in the B\&C Search Tree (IEEE 33-Node Test System)}
  \vspace{-6pt}
  \setlength{\tabcolsep}{7.5pt}
    \begin{tabular}{|c|c|c|c|c|}
      \hline
    \multirow{2}[0]{*}{Avg. $T_s$} &  \multirow{2}[0]{*}{Avg. $T_m$} &
    \multicolumn{3}{c|}{Number of cases with:} \\
    \cline{3-5}
  & &  $T_s< T_m$  & $T_s= T_m$  & $T_s > T_m$\\
      \hline
  0.49 s & 0.56 s & 8911/10000 & 0/10000 & 1089/10000\\
  \hline
  \hline
      \multirow{2}[0]{*}{Avg. $N_s$} &  \multirow{2}[0]{*}{Avg. $N_m$} &
      \multicolumn{3}{c|}{Number of cases with:} \\
      \cline{3-5}
      & &  $N_s< N_m$  & $N_s= N_m$  & $N_s > N_m$\\
        \hline
      407 & 143 & 402/10000 & 21/10000 & 9577/10000  \\
      \hline
\multicolumn{5}{l}{\tabincell{l}{$T_s/T_m$ ({\color{black}resp. }$N_s/N_m$):
Computation time ({\color{black}resp. }the number of ex-\\plored nodes
in
the B\&C search tree) of the proposed MG formation \\ model using the
single-commodity flow-based method$/$directed multi-\\commodity flow-based
method to formulate constraint~\eqref{SF1}.}}
    \end{tabular}%
  \label{node_time}%
  \vspace{-8pt}
\end{table}%

\subsection{IEEE 123-Node Test System \cite{123system}}

\begin{figure}[t!] 
  \centering
  \includegraphics[width=3.49in]{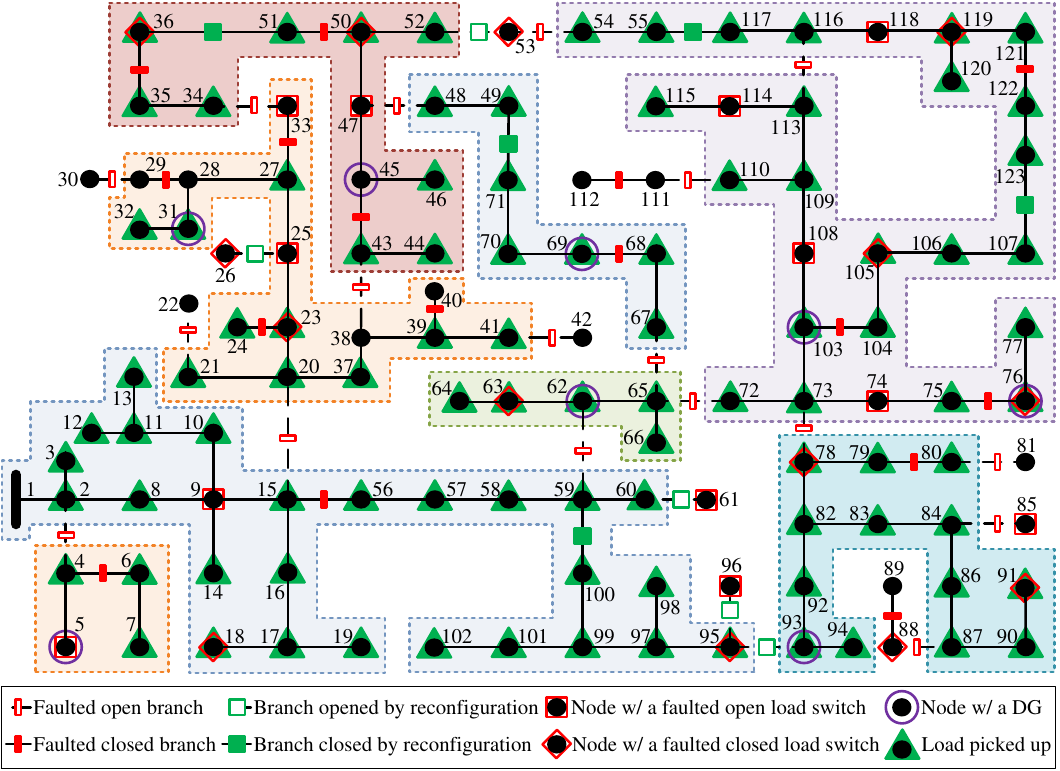}
  \vspace{-18pt}
  \caption{An illustrative case of MG formation on IEEE 123-node test system~(use our proposed model).}
  \label{123node_base_case}
\end{figure}

\begin{table}[t!]
  \centering
  \caption{Computation Time and the Number of Infeasible Cases (IEEE 123-Node Test System)}
  \vspace{-6pt}
    \setlength{\tabcolsep}{5.8pt}
    \begin{tabular}{|c|c|c|c|}
    \hline
    --------- & Model in \cite{Chen2016} & Model in \cite{Ding2017} & Proposed model \\
    \hline
    Avg. comput. time & 13.41 s & 4.36 s & 1.79 s \\
    \hline
    Infeasible cases &
    109/10000 &
    196/10000  &
    2/10000 \\
    \hline
    \end{tabular}%
  \label{infeasible_number_2}%
\end{table}%

\begin{table}[t!]
  \centering
  \caption{Summary Statistics for the Restored Loads (IEEE 123-Node Test System)}
  \vspace{-6pt}
  \setlength{\tabcolsep}{2pt}
    \begin{tabular}{|c|c|c|c|c|c|c|}
    \hline
    \multirow{2}[0]{*}{------} & \multicolumn{5}{c|}{Restored loads (kW)}      & Box plots \\
    \cline{2-7}
     & avg.  & std.  & max.  & med. & min.  & \multirow{4}[0]{*}{\raisebox{-0.1in}{\includegraphics[width=0.94in]{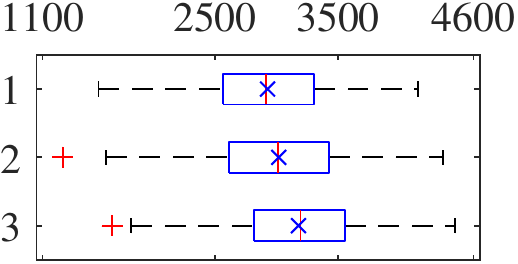}}} \\
     \cline{1-6}
    \tabincell{c}{1: Model in \cite{Chen2016}}     &  2917  &  482  &    4150   &   2915    &   1555    &  \\
    \cline{1-6}
    \tabincell{c}{2: Model in \cite{Ding2017}}     &  3008 &  524  &   4350    &    3010   &    1265   &  \\
    \cline{1-6}
    \tabincell{c}{3: Proposed model}     &   3189    &   494   & 4450 &   3195   &   1665    &  \\
    \hline
    \end{tabular}%
  \label{box_obj_2}%
\end{table}%

\begin{table}[t!]
  \centering
  \caption{Summary Statistics for the DG Capacity Utilization Rate (IEEE 123-Node Test System)}
\vspace{-6pt}
  \setlength{\tabcolsep}{0.99pt}
    \begin{tabular}{|c|c|c|c|c|c|c|}
    \hline
    \multirow{2}[0]{*}{------} & \multicolumn{5}{c|}{DG capacity utilization rate}      & Box plots \\
    \cline{2-7}
     & avg.  & std.  & max.  & med. & min.  & \multirow{4}[0]{*}{\raisebox{-0.1in}{\includegraphics[width=1.01in]{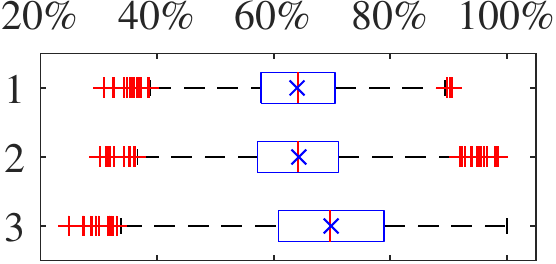}}} \\
     \cline{1-6}
    \tabincell{c}{1: Model in \cite{Chen2016}}     &  63.0\%  &  9.1\%  &    90.5\%  &   63.1\%    &   30.8\%    &  \\
    \cline{1-6}
    \tabincell{c}{2: Model in \cite{Ding2017}}     &  63.3\% &  10.6\%  &   97.4\%    &    63.1\%   &    30.1\%   &  \\
    \cline{1-6}
    \tabincell{c}{3: Proposed model}     &   70.9\%    &   12.4\%   & 100\% &   70.7\%   &   24.9\%    &  \\
    \hline
    \end{tabular}%
  \label{box_DG_2}%
\end{table}%

\begin{figure}[t!]
\hspace*{\fill}%
  \begin{minipage}[t!]{4.4cm}
    \centering
   \vspace{0pt}
    \includegraphics[height=1.35in]{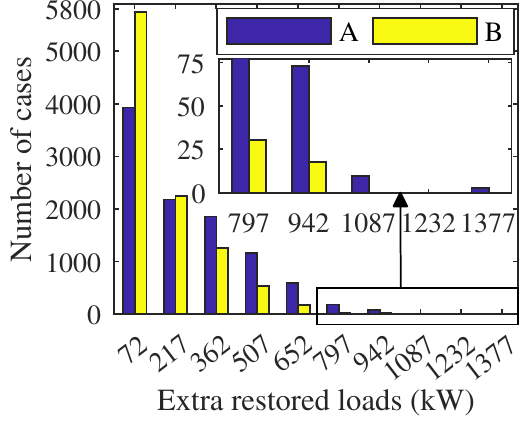}
    \label{extra_loads_2}
  \end{minipage}%
\hfill
  \begin{minipage}[t!]{4.4cm}
    \centering
    \vspace{0pt}
    \includegraphics[height=1.35in]{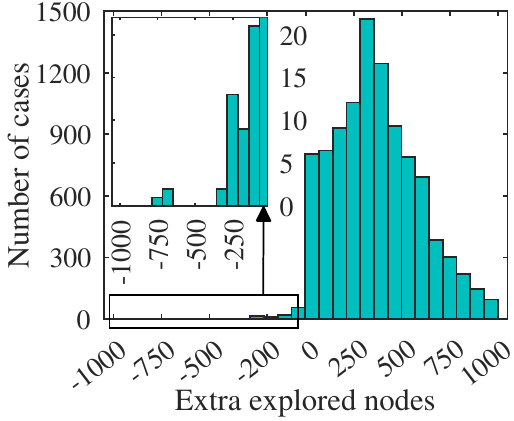}
    \label{extra_nodes_2}
  \end{minipage}
  \hspace*{\fill}
{\hfill}
    \vspace{-3pt}
  \caption{Left: Histograms of the extra restored loads of our proposed~MG~formation
  model, compared with the model in \cite{Chen2016} (A),~and~compared~with~the~model in \cite{Ding2017} (B).
  Right: Histogram of the extra~explored nodes in the B\&C search tree of
  the proposed MG formation model~using~the~single-commodity flow-based method to formulate
  constraint \eqref{SF1}, compared~with~that~using~the~directed~multicommodity
  flow-based method. (IEEE 123-node~test~system.)}
  \label{histograms_2}
\end{figure}

Fig. \ref{123node_base_case} provides an illustrative case
based on a 60-fault~sce-

\noindent
nario
using the proposed MG formation model.
It forms 7 MGs powered by 8 DGs, and a sub-grid powered by the substation.
The DG capacity utilization rate is 75.0\%,
and 3730kW loads are restored.
We again run 10000 cases on
random scenarios of DS faults.
Table \ref{infeasible_number_2}-\ref{node_time_2}
and Fig. \ref{histograms_2} show results
similar to the previous system.
That is, the proposed radiality constraints
and MG formation model's superiority is also established
on this larger system.
For example, our model matches DGs and loads in a more
coordinated manner, thus
achieving better service restoration.
For space limit, we do not go into the details.
In short, our MG formation
model has extended feasibility and enhanced optimality due
to the proposed radiality constraints' fully enabling topological
and some related flexibilities of DSs.

It is worth mentioning that, as in
Table \ref{node_time_2}, $T_m$ is averagely shorter
than $T_s$ for this larger system.
Specifically, in the~cases with $T_m<T_s$, $T_m$
and $T_s$ are 1.26 s and 5.59 s on~average, respectively; $N_m$ and
$N_s$ are 89 and 937 on average,~respectively.
That is, the computation time~of~the~tighter~version~of~our~MG formation model is not only shorter on average,
but also~more consistent.
It greatly reduces the computation time especially for cases requiring
the less tight model to explore much~more B\&C
search tree nodes. Thus, the tighter version~of~the~proposed radiality constraints
can allow both more~flexible DS operation and more efficient computation for
reconfiguration-related
optimization problems, e.g., MG formation here.

\begin{table}[t!]
  \centering
  \caption{Computation Time and the Number of Explored Nodes in the B\&C Search Tree (IEEE 123-Node Test System)}
  \setlength{\tabcolsep}{7.5pt}
    \begin{tabular}{|c|c|c|c|c|}
    \hline
  \multirow{2}[0]{*}{Avg. $T_s$} &  \multirow{2}[0]{*}{Avg. $T_m$} &
  \multicolumn{3}{c|}{Number of cases with:} \\
  \cline{3-5}
& &  $T_s< T_m$  & $T_s= T_m$  & $T_s > T_m$\\
    \hline
1.79 s & 1.54 s & 7348/10000 & 0/10000 & 2652/10000\\
\hline
\hline
\multirow{2}[0]{*}{Avg. $N_s$} &  \multirow{2}[0]{*}{Avg. $N_m$} &
\multicolumn{3}{c|}{Number of cases with:} \\
\cline{3-5}
& &  $N_s< N_m$  & $N_s= N_m$  & $N_s > N_m$\\
  \hline
496 & 71 & 108/10000 & 134/10000 & 9758/10000  \\
\hline
\multicolumn{5}{l}{\tabincell{l}{$T_s/T_m$ ({\color{black}resp. }$N_s/N_m$):
Computation time ({\color{black}resp. }the number of ex-\\plored nodes
in
the B\&C search tree) of the proposed MG formation \\ model using the
single-commodity flow-based method$/$directed multi-\\commodity flow-based
method to formulate constraint~\eqref{SF1}.}}
    \end{tabular}%
  \label{node_time_2}%
        \vspace{-3pt}
\end{table}%

\section{Conclusion}\label{conclusions_final}

This work proposes a new method for formulating
radiality constraints that fully enables
topological and some other related flexibilities
in DS reconfiguration optimization
problems.
It is specifically applied to resilient post-disaster MG formation
to attain
extended feasibility and enhanced optimality.
As verified in case studies, compared to the
existing MG formation models,
our model based on the proposed radiality constraints
achieves a higher resilience enhancement via
more coordinated utilization of DS flexibilities, etc.,
and reduces the computational~complexity.
Future work includes
exploring the effects
of increased topological flexibilities
on other DS performance metrics, etc.

\ifCLASSOPTIONcaptionsoff
  \newpage
\fi



%
\bibliographystyle{IEEEtran}
\bibliography{IEEEabrv,references}

%








\end{document}